\documentclass[a4, 12pt]{amsart}
\usepackage{amssymb}
\usepackage{amstext}
\usepackage{amsmath}
\usepackage{amscd}
\usepackage{latexsym}
\usepackage{amsfonts}
\usepackage[all]{xy}

\theoremstyle{plain}
\newtheorem{thm}{Theorem}[section]
\newtheorem*{thm*}{Theorem}
\newtheorem*{cor*}{Corollary}

\newtheorem{prop}[thm]{Proposition}
\newtheorem{lem}[thm]{Lemma}
\newtheorem{cor}[thm]{Corollary}
\newtheorem{claim}{Claim}
\newtheorem*{claim*}{Claim}

\theoremstyle{definition}

\newtheorem{ex}[thm]{Example}

\newtheorem{conj}[thm]{Conjecture}

\theoremstyle{remark}
\newtheorem*{pf}{{\sl Proof}}
\newtheorem*{tpf}{{\sl Proof of Theorem 1.1}}
\newtheorem*{cpf1}{{\sl Proof of Claim 1}}
\newtheorem*{cpf2}{{\sl Proof of Claim 2}}
\newtheorem*{cpf}{{\sl Proof of Claim}}

\numberwithin{equation}{thm}
\def\Hom{\mathrm{Hom}}

\def\Ext{\mathrm{Ext}}

\def\mod{\mathrm{mod}}

\def\e{\mathrm{e}}
\def\m{\mathfrak m}
\def\n{\mathfrak n}

\def\Gr{\mathrm G}
\def\FC{\mathrm F}

\newcommand{\rma}{\mathrm{a}}

\newcommand{\rme}{\mathrm{e}}

\newcommand{\rmr}{\mathrm{r}}

\newcommand{\rmv}{\mathrm{v}}

\newcommand{\rmF}{\mathrm{F}}
\newcommand{\rmG}{\mathrm{G}}
\newcommand{\rmH}{\mathrm{H}}
\newcommand{\rmI}{\mathrm{I}}

\newcommand{\calR}{\mathcal{R}}

\newcommand{\fkm}{\mathfrak{m}}

\def\R{{\mathcal R}}

\begin{document}

\setlength{\baselineskip}{20pt}

\title{Quasi-socle ideals in  a Gorenstein local ring}
\author{Shiro Goto}
\address{Department of Mathematics, School of Science and Technology, Meiji University, 1-1-1 Higashimita, Tama-ku, Kawasaki 214-8571, Japan}
\email{goto@math.meiji.ac.jp}
\author{Naoyuki Matsuoka}
\address{Department of Mathematics, School of Science and Technology, Meiji University, 1-1-1 Higashimita, Tama-ku, Kawasaki 214-8571, Japan}
\email{matsuoka@math.meiji.ac.jp}
\author{Ryo Takahashi}
\address{Department of Mathematical Sciences, Faculty of Science, Shinshu University, 3-1-1 Asahi, Matsumoto, Nagano 390-8621, Japan}
\email{takahasi@math.shinshu-u.ac.jp}
\thanks{
\endgraf
}
\maketitle
\begin{abstract}
This paper explores the structure of quasi-socle ideals $I=Q:\m^2$ in a Gorenstein local ring $A$, where $Q$ is a parameter ideal and $\m$ is the maximal ideal in $A$. The purpose is to answer the problems of when $Q$ is a reduction of $I$ and when the associated graded ring $\rmG(I) = \bigoplus_{n \geq 0}I^n/I^{n+1}$ is Cohen-Macaulay. Wild examples are explored.
\end{abstract}

\section{Introduction}

The purpose of this paper is to prove the following theorem.

\begin{thm}\label{gmt}
Let $(A,\m)$ be a Gorenstein local ring with $\dim A > 0$ and assume that $\rme_\m^0(A) \geq 3$, where $\rme_\m^0(A)$ denotes the multiplicity of $A$ with respect to the maximal ideal $\m$. Then for every parameter ideal $Q$ in $A$, one has the following, where $I = Q:\m^2$.
\begin{enumerate}
\item[{\rm (1)}] $\m ^2I=\m ^2Q$ and $I^3 =QI^2$.
\item[{\rm (2)}] The associated graded ring $\Gr(I)$ of $I$ and the fiber cone $\rmF(I)$ of $I$ are both Cohen-Macaulay rings.
\end{enumerate}
Hence, the Rees algebra $\calR(I)$ of $I$ is also a Cohen-Macaulay ring, if $\dim A \geq 3$.
\end{thm}

\noindent
Here we define\\
\hspace{10mm}$\calR(I) = A[IT] \ \ \subseteq A[T], \\
\hspace{10mm}\calR'(I) =A[IT, T^{-1}] \ \ \subseteq A[T, T^{-1}],\\
\hspace{10mm}\rmG(I) = \calR'(I)/T^{-1}\calR'(I), ~\textup{and}\\
\hspace{10mm}\rmF(I) = \calR(I)/\m \calR(I) \ \ (\cong \rmG(I)/\m \rmG(I))$ \\
with $T$ an indeterminate over $A$.

Our Theorem 1.1 is a generalization  of the following result of  A. Corso, C. Polini, C. Huneke, W. V. Vasconcelos, and the first author.

\begin{thm}[\cite{CHV, CP1, CP2, CPV, G}]\label{q=1}
Let $(A,\m)$ be a Cohen-Macaulay local ring with $d = \dim A > 0$.
Let $Q$ be a parameter ideal in $A$ and let $I=Q:\m$.
Then the following three conditions are equivalent to each other.
\begin{enumerate}
\item[{\rm (1)}] $I^2 \ne QI$.
\item[{\rm (2)}] $\overline{Q} = Q$, that is the parameter ideal $Q$ is integrally closed in $A$.
\item[{\rm (3)}] $A$ is a regular local ring which contains a regular system $x_1, x_2, \cdots , x_d$ of parameters such that $Q = (x_1, \cdots , x_{d-1}, x_d^q)$ for some integer $q > 0$.
\end{enumerate}
Consequently, if $(A,\m)$ is a Cohen-Macaulay local ring which is not  regular, then $I^2=QI$ for every parameter ideal $Q$ in $A$, so that $\rmG(I)$ and $\rmF(I)$ are both Cohen-Macaulay rings, where $I=Q:\fkm$. The Rees algebra $\calR (I)$ is also a Cohen-Macaulay ring, if $d = \dim A \geq 2$.
\end{thm}

The present research aims at a natural generalization of Theorem 1.2 but here we would like to note  that there might be other directions of generalization.  In fact, the equality $I^2=QI$ in Theorem 1.2 remains true in certain cases, even though the base local rings $A$ are not Cohen-Macaulay. For example, the first author and H. Sakurai investigated the case where $A$ is a Buchsbaum local ring and gave the following. See \cite{GSa1, GSa3} for further developments of this direction.

\begin{thm}[\cite{GSa2}, cf. \cite{GN}]\label{1.5}
Let $(A,\fkm)$ be a Buchsbaum local ring and assume that either $\dim A \ge 2$ or $\dim A=1$ but $\rme_\fkm^0(A) \ge 2$. Then there exists an integer $n >0$ such that for every parameter ideal $Q$ of $A$ which is contained in $\fkm^n$, one has the equality $I^2=QI$, so that the graded rings $\rmG(I)$ and $\rmF(I)$ are Buchsbaum rings, where $I=Q:\fkm$.
\end{thm}

The researches \cite{CHV, CP1, CP2, CPV} originate at the study of linkage of ideals. If $A$ is a Cohen-Macaulay local ring and $I$ is an equimultiple Cohen-Macaulay ideal in $A$ of reduction number one, then the associated graded ring $\mathrm{G}(I)$ is Cohen-Macaulay and, so is the Rees algebra $\R(I)$, provided $\operatorname{ht}_AI \geq 2$. One knows the number and degrees of the defining equations of $\R(I)$ also, so that one can understand fairly explicitly the process of desingularization of $\operatorname{Spec} A$ along the subscheme $\mathrm{V}(I)$. This motivated the ingenious  research of C. Polini and B. Ulrich \cite{PU}, where they posed, with many other important results, the following conjecture

\begin{conj}[\cite{PU}]\label{conj}
Let $(A, \m)$ be a Cohen-Macaulay local ring with $\dim A \geq 2$. Assume that $\dim A \geq 3$ when  $A$ is regular. Let $q \geq 2$ be an integer and  $Q$  a parameter ideal in $A$ such that $Q \subseteq \m^q$. Then
 $$Q:\m^q \subseteq \m^q.$$ 
\end{conj}

\noindent
and H.-J. Wang \cite{W} recently  settled  this conjecture in the following way.

\begin{thm}[\cite{W}]\label{w}
Let $(A, \m)$ be a Cohen-Macaulay local ring with $d=\dim A \geq 2$. Let $q \geq 1$ be an integer and $Q$  a parameter ideal in $A$. Assume that $Q \subseteq \m^q$ and put $I = Q : \m^q$. Then $$I \subseteq \m^q, ~~\m^q I = \m^q Q,~~~~\operatorname{and}~~~~I^2 = QI,$$ provided that $A$ is not regular, if $d \geq 2$ and that $q \geq 2$, if $d \geq 3$.
 \end{thm}

Wang's result Theorem \ref{w} is certainly closely related to our Theorem \ref{gmt}, although  Theorem \ref{w}, apparently, does not cover our Theorem \ref{gmt}. The two researches were performed independently and our proof of method is,  heavily depending on the facts that the base ring $A$ is Gorenstein, $q = 2$, and $\e_{\m}^0(A) \geq 3$, totally different from Wang's method, and despite the restrictions, our Theorem \ref{gmt} holds true for every parameter ideal $Q$ in $A$, even in the case where $\dim A = 1$. For this reason, Theorem \ref {gmt} may have its own significance, suggesting a possible modification of the Polini-Ulrich conjecture.


We now explain  how this paper is organized. Section 2 is devoted to some preliminary steps, which we will need later to prove Theorem 1.1.
Theorem 1.1 will be proven in Section 3. Our method of proof is, unfortunately, applicable only to the case where the local ring $A$ is Gorenstein and the situation seems totally different, unless $A$ is Gorenstein. In order to show that the non-Gorenstein case of dimension 1 is rather wild, we shall explore three examples in the last Section 4. 
One of them  will show the quasi-socle ideals $I=Q:\m^2$ are never integral over parameter ideals $Q$ in certain Cohen-Macaulay local rings $A$ of dimension 1, even though $\rme_\m^0(A) \geq 2$. The other two will show that  unless $A$ is a Gorenstein ring, one can not expect that $\rmr_Q(I) \leq 2$, even if $I$ is integral over $Q$, where $$\rmr_Q(I) = \min \{n \geq 0 \mid I^{n+1}=QI^n \}$$ 
denotes the reduction number of the ideal $I=Q:\m^2$ with respect to $Q$.

Unless otherwise specified, in what follows, let $(A,\m)$ be a Gorenstein local ring with $\dim A=d$. We denote by $\rme_\m^0(A)$ the multiplicity of $A$ with respect to the maximal ideal $\m$. Let $Q =(a_1, a_2, \cdots, a_d)$ be a parameter ideal in $A$ generated by the system $a_1, a_2, \cdots, a_d$ of parameters in $A$. For each finitely generated $A$-module $M$ we denote by $\mu_A(M)$ and $\ell_A(M)$, respectively,  the number of elements in a minimal system of generators for $M$ and the length of $M$. Let $\rmv(A) = \ell_A(\m/\m^2)$ stand for  the embedding dimension of $A$. 

\section{Preliminaries}

Let $A$ be a Gorenstein local ring with the maximal ideal $\m$.  The purpose of this section is to summarize some preliminaries, which we need in Section 3 to prove Theorem 1.1. Let us  begin with the case where $\dim A = 0$.

Suppose that $\dim A =0$. Let $n=\rmv(A)> 0$ and let $x_1,x_2, \cdots, x_n$ be a system of generators for $\m$. 
We choose  a socle element $z$ in $A$. Hence $0 \ne z \in \m$ and $\m z =(0)$. Let $I = (0) : \m^2$. We then  have the following. 

\begin{lem}\label{d=0}
There exist elements $y_1,y_2,\cdots, y_n \in A$ such that $x_iy_j=\delta _{ij}z$ for all integers $1\le i,j\le n$. We furthermore have the following.
\begin{enumerate}
\item[{\rm (1)}]
$I =(y_1, y_2, \cdots, y_n)$, $\mu (I)=n$, and $\ell_A(I) = n+1$.
\item[{\rm (2)}]
If $n > 1$, then $I \subsetneq A$.
\end{enumerate}
\end{lem}

\begin{pf}
The existence of elements $y_1, y_2, \cdots, y_n$ is exactly  the dual basis lemma. Let us note a brief proof for the sake of completeness. Let $1\le j\le n$ be an integer. We look at the following diagram 
\[\xymatrix{
\m \ar[rrr]^{\iota} \ar[d]^\varepsilon & & & A \ar@{-->}[d]^{f = \widehat{y_j}}  \\
\m /\m^2 \ar[r]^p & A/\m \ar[r]^{\overset{h}{\sim}} & (z) \ar[r]^\iota & A
}\]
of $A$-modules, where $\varepsilon$ is the canonical epimorphism, $p$ is the projection map such that $p(\overline{x_i}) =\delta_{ij}$ for all $1 \leq i \leq n$ where $\overline{x_i} = x_i~\mod~\m^2$ denotes the image of $x_i$ in $\m/ \m^2$ and $\delta _{ij}$ is Kronecker's delta, $h$ is the isomorphism of vector spaces over $A/\m$ defined by $h(1) = z$, and $\iota$'s denote the embedding maps.
Then, since the ring $A$ is self-injective, we have a homothety map $f = \widehat{y_j} :A\to A$ with $y_j \in A$ such that the above diagram is commutative. Hence $x_iy_j =\delta _{ij}z$ for all integers $1 \leq i, j \leq n$.
We put $J=(y_1,y_2,\cdots,y_n)$. Then $J \subseteq I=(0) : \m^2$, because $\m z=(0)$ and $x_iy_j = \delta_{ij}z$. We have $\ell_A(I) = n+1$, since $$I \cong \Hom _A(A/\m^2,A)~\ ~\textup{and}~\ ~\ell_A(A/\m^2)=n+1.$$ Therefore, to see that $I=J$, we have only to show $\ell_A(J) = n+1$, or equivalently $\ell_A(J/(z)) = n$.
Let $\{b_j\}_{1 \le j \le n}$ be elements in $A$ and assume that $\sum_{j=1}^nb_jy_j \in (z)$.
Then $$b_iz = b_i(x_iy_i) = x_i{\cdot}\sum_{j=1}^nb_jy_j= 0.$$ Hence $b_i \in \m$.
Thus the images of $\{y_j\}_{1\le j \le n}$ in $J/(z)$ form a basis of the vector space $J/(z)$ over $A/\m$, so that $\mu_A(J/(z)) = \ell_A(J/(z))=n$. Hence $\ell_A(J) = n + 1$ and assertion (1) follows. Assertion (2) is now obvious.
\qed
\end{pf}

For the rest of this section we throughout assume that $d = \dim A > 0$. Let $Q=(a_1, a_2, \cdots, a_d)$ be a parameter ideal in $A$ generated by a system $a_1, a_2, \cdots, a_d$ of parameters for $A$ and let $I = Q:\m^2$. We assume $n=\rmv(A/Q) > 0$ and  write $\m = Q + (x_1, x_2, \cdots, x_n)$ with $x_i \in A$. Then $\m I \subseteq Q : \m$ and $\m I \not\subseteq Q$ (recall that $Q \ne \m$, since $n> 0$). Let us choose $z \in \m I$ so that $z \not\in Q$, whence $$Q:\fkm = Q + \fkm I = Q+(z).$$ Then, applying Lemma \ref{d=0} to the Artinian local ring $A/Q$, we get the elements $y_1, y_2, \cdots , y_n \in A$ such that 
$x_iy_j \equiv \delta_{ij}z ~~\mod ~~Q$
for all integers $1 \le i, j \le n$.
Hence $$I = Q + (y_1, y_2, \cdots , y_n ), ~~~\ \mu_A(I/Q) = n, ~~~\ \textup{and}~~~\ \ell_A(I/Q) = n+1,$$ 
so that we have $\mu_A(I) \leq n+d$. 

We now look at the following inclusions

\[\xymatrix{
I \ar@{-}[dd]_{n+1} \ar@{-}[rd]^{\mu_A(I) \le n+d} &  \\
  & \fkm I \ar@{-}[dd]^1\\
Q \ar@{-}[dd]_d &  \\
  & \fkm I \cap Q \ar@{-}[ld]\\
\fkm Q &
}\]

\vspace{5mm}
\noindent
and notice that $[Q + \fkm I]/Q \cong \fkm I/[\fkm I \cap Q]$.
Then $\ell_A(\fkm I/[\fkm I \cap Q]) = 1$ since $Q:\fkm = Q + \fkm I$, so that we have
$$
\mu_A(I) = n+d ~\Longleftrightarrow ~\fkm I \cap Q = \fkm Q.
$$ 

We furthermore have the following.

\begin{prop}\label{4}
Suppose that $n=\rmv(A/Q) > 1$.
Then the following four conditions are equivalent to each other.
\begin{enumerate}
\item[{\rm (1)}]
$I \subseteq \overline{Q}$.
\item[{\rm (2)}]
$\m I\cap Q=\m Q$.
\item[{\rm (3)}]
$\mu _A(I)=n+d$.
\item[{\rm (4)}]
$\m ^2I=\m ^2Q$.
\end{enumerate}
Here $\overline{Q}$ denotes the integral closure of $Q$.
\end{prop}

\begin{pf}
The implication (1) $\Rightarrow$ (2) is clear, since $Q$ is a minimal reduction of $I$. The equivalence (2) $\Longleftrightarrow$ (3) follows from the above observation. 

(4) $\Rightarrow$ (1) This is well-known (cf. \cite{NR}). Use the determinantal trick. 

(2) $\Rightarrow$ (4) Because $z \in \fkm I \subseteq Q:\fkm = Q+(z)$, we get
\begin{eqnarray*}
\fkm I &=& (\fkm I \cap Q) + (z)\\
&=& \fkm Q + (z).
\end{eqnarray*}
Therefore, in order to see the equality  $\fkm^2 I = \fkm^2 Q$, we have only to show that $$\fkm z \subseteq \fkm^2 Q.$$
Since $z \in \fkm I \subseteq \fkm^2$ (recall that $I \ne A$; cf. Lemma \ref{d=0} (2)), we get $Qz \subseteq \frak{m}^2Q$. Hence, because  $\fkm = Q+(x_1, x_2, \cdots ,x_n)$, it suffices to show that $x_{\ell}z \in \frak{m}^2Q$ for every $1 \leq \ell \leq n$. Choose an integer $1 \leq i \leq n$ so that $i\ne \ell$ and write $z=x_iy_i+q_i$ with $q_i \in Q$.
Then $x_{\ell}z = x_i(x_\ell y_i) + x_\ell q_i$. Because $q_i = z - x_iy_i \in \fkm I \cap Q = \fkm Q$ and $x_\ell y_i \in \fkm I \cap Q =  \fkm Q$, we certainly have $x_\ell z \in \fkm^2 Q$. Thus  $\fkm^2 I = \fkm^2 Q$.
\qed
\end{pf}

As a consequence of Proposition \ref{4} we have the following.

\begin{cor}\label{qi}
Suppose that $n=\rmv(A/Q) > 1$ and that $I$ is integral over $Q$. Then 
\begin{enumerate}
\item[{\rm (1)}]
$Q^i\cap I^{i+1}=Q^iI$ for all integers $i \geq 1$. Hence $I^2=QI$ if $I \subseteq \m^2$.
\item[{\rm (2)}]
$(a_1)\cap I^2=a_1I$.
\item[{\rm (3)}]
$I^2=QI$ if $Q \subseteq \m^2$.
\end{enumerate}
\end{cor}

\begin{pf}
(1) The second assertion follows from the first, since $I^2 \subseteq \m^2 I \subseteq Q$. To see the first assertion, notice that $\m^2 I^{i+1} = \m^2 Q^{i+1}$, since $\m^2 I = \m^2 Q$ by Proposition \ref{4}. Let $f \in Q^i\cap I^{i+1}$ and write 
$$f =\sum _{i_1+i_2+\cdots +i_d=i}a_1^{i_1}a_2^{i_2}\cdots a_d^{i_d}f_{i_1i_2\cdots i_d}$$
with $f_{i_1i_2\cdots i_d}\in A$.
Let $\alpha \in \m^2$. We then have $$\alpha f = \sum _{i_1+i_2+\cdots +i_d=i}a_1^{i_1}a_2^{i_2}\cdots a_d^{i_d}(\alpha f_{i_1i_2\cdots i_d}) \in \m^2 I^{i+1} \subseteq Q^{i+1}.$$ Hence $\alpha f_{i_1i_2\cdots i_d} \in Q$ because $a_1, a_2, \cdots ,a_d$ is an $A$-regular sequence, so that $f_{i_1i_2\cdots i_d} \in I$. Thus $f \in Q^iI$, whence $Q^i\cap I^{i+1}=Q^iI$.

(2) Let $f \in (a_1)\cap I^2$ and write $f =a_1g$ with $g \in A$.
Then for all $\alpha \in \m^2$, we have $\alpha f = a_1(\alpha g) \in \m^2 I^2 \subseteq Q^2$. Hence $\alpha g \in Q$ so that $g\in I$, and so $f \in a_1 I$. Thus $(a_1)\cap I^2=a_1I$.

(3) Let us prove the assertion by induction on $d$.
Assume that $d=1$. Let $b\in \m^2$ be a non-zerodivisor in $A$. Then, thanks to the isomorphisms 
$$[(b):\m^2]/(b) \cong \Hom _A(A/\m^2,A/(b)) \cong\Ext _A^1(A/\m ^2,A)$$
of $A$-modules,
\noindent
we see the length $\ell_A([(b):\m^2]/(b)) = \ell_A(\Ext _A^1(A/\m ^2,A))$ is independent of the choice of the element $b \in \m^2$.
We put $a = a_1$. Let $Q'=(a^2)$ and $I'=Q':\m^2$.
Let $$\varphi : A/(a) \to A/(a^2)$$ be the monomorphism defined by $\varphi(\overline x)=\overline{ax}$, where $\overline{*}$ denote the images of the corresponding elements $x$ and $ax$. Then $\varphi(I/(a)) = I'/(a^2)$, since $\varphi(I/(a)) \subseteq I'/(a^2)$ and $\ell_A(I/(a))=\ell_A(I'/(a^2))$ (recall that $a \in \m^2$). 
Therefore $$(\sharp) \ \ \ \ \ \ I'=aI+(a^2)=aI,$$ whence $\mu _A(I')=\mu _A(I)=n+1$, where the last equality follows from Proposition \ref{4}. Hence $I'$ is also  integral over $Q'$ by Proposition \ref{4}, because $\rmv(A/Q')=\rmv(A) = \rmv(A/Q)=n > 1$.
Therefore $(I')^2=a^2I'$ by assertion (1), since $I' \subseteq \m^2$. Hence by equality $(\sharp)$ we get  $a^2I^2 = (I')^2 = a^2I' = a^3I$, so that $I^2 = aI$.

Assume now that $d\ge 2$ and that our assertion holds true for $d-1$.
Let $\overline A=A/(a_1)$, $\overline\m = \m /(a_1)$, $\overline Q=Q/(a_1)$, and $\overline I=I/(a_1)$.
Then $\overline Q: \overline \m^2 =\overline I$, $\rmv(\overline A/\overline Q)=\rmv(A/Q) = n > 1$, and $\overline{I}$ is integral over $\overline{Q}$. Hence the hypothesis of induction on $d$ yields that $\overline I^2=\overline Q~\overline I$, since $\overline Q \subseteq \overline\m ^2$. Thus $I^2 \subseteq QI+(a_1)$.
Therefore 
$$I^2=[QI+(a_1)]\cap I^2=QI+[(a_1)\cap I^2]=QI+a_1I=QI$$
 by assertion (2).
\qed
\end{pf}

\begin{cor}\label{5}
Suppose $\rmv(A/Q) > 1$ and $I$ is integral over $Q$. Then $I\subseteq\m^2$ if $Q\subseteq\m^2$.
\end{cor}

\begin{pf}
Suppose $Q\subseteq\m^2$. Then $I^2 \subseteq Q$ since $I^2=QI$ by Corollary \ref{qi} (3). On the other hand we have $Q:(Q:\m^2)=\m^2$, because $Q$ is a parameter ideal in the Gorenstein local ring $A$. Hence $I \subseteq Q:I = Q:(Q:\m^2) = \m^2$ as is claimed.
\qed
\end{pf}

Unless $Q \subseteq \fkm^2$, the equality $I^2=QI$ does not necessarily hold true. Let us note one example.

\begin{ex}
Let $H=\left<6, 7, 15 \right>$ be the numerical semi-group generated by $6, 7, 15$ and let $A=k[[t^6, t^7, t^{15}]]~\subseteq~k[[t]]$, where $k[[t]]$ denotes the formal power series ring with one indeterminate $t$ over a field $k$ . Then $A$ is a Gorenstein local ring with $\dim A = 1$. Let $0 < s \in H=\left<6, 7, 15\right>$, $Q=(t^s)$ in $A$, and $I=Q:\fkm^2$. Then $I$ is integral over $Q$ and $\rmr_Q(I) \le 2$. However,  $I^2=QI$ if and only if $s \ne 7.$
\end{ex}

\begin{pf}
Let $n \in H$. Then it is direct to check that $t^n \in I$ if and only if $n =s, s+6, s+7, s+8,~\textup{or}~s+\ell~\textup{for some}~ 12 \le \ell \in \Bbb Z$. Thanks to this observation, we get $I = (t^s, t^{s + 8}, t^{s+16}, t^{s+17})$ if $s \ge 12$ but $s \ne 15$. We also have $I = (t^6, t^{14}, t^{22})$ if $ s = 6$, $I =(t^7, t^{15}, t^{24})$ if $s = 7$, and $I = (t^{15}, t^{31}, t^{32})$ if $s = 15$. Hence $I \subseteq t^sk[[t]] \cap A$,  so that $I$ is integral over $Q=(t^s)$,  in any case. It is routine to check that $I^2 = QI$ when $s \ne 7$. If $s = 7$, then $I^3 = QI^2$ but $I^2= QI +(t^{30})$ and $\ell_A(I^2/QI) = 1$, whence $I^2 \ne QI$.
\qed
\end{pf}

Here let us note one example to clarify our arguments.

\begin{ex}
Let $(A, \m)$ be a regular local ring with $d = \dim A \ge 2$ 
and let $x_1,x_2, \cdots,x_d$ be a regular system of parameters of $A$.
Let $c_i \ge 2$ $(1 \le i\le d)$ be integers and put $Q=(x_1^{c_1}, x_2^{c_2},\cdots,x_d^{c_d})$. Let  $I=Q:\m^2$. We then have the following.
\begin{enumerate}
\item[(1)]
The following conditions are equivalent.
\begin{enumerate}
\item[(i)]
$I \not\subseteq \overline{Q}$.
\item[(ii)]
$d=2$ and $\min\{c_1, c_2\} = 2$.
\end{enumerate}
\item[(2)]
$I^2=QI$ if $I \subseteq \overline{Q}$.
\end{enumerate}
Here $\overline{Q}$ denotes the integral closure of $Q$.
\end{ex}
\begin{pf}
Let  $z=\prod_{i=1}^d x_i^{c_i-1}$, $a_i=x_i^{c_i}$, and $y_i=\displaystyle\frac{z}{x_i}$ for each $1\le i\le d$.
Then $Q:\m = Q+ (z)$ and $x_iy_j \equiv\delta_{ij}z$ modulo $Q$ for all integers $1\le i,j\le d$.
Hence $I=Q+(y_1, y_2, \cdots,y_d)$ and $\mu_A(I/Q) = d$ by Lemma \ref{d=0}. We put $J=(y_1, y_2, \cdots,y_d)$. 

Suppose now that $I\not\subseteq \overline Q$.
Then, since $\rmv(A/Q) = d > 1$, by Proposition \ref{4} we have $\mu_A(I)<2d$. Hence $a_i \in (a_j \mid 1 \le j \le d, j \ne i) + J$ for some $1 \le i \le d$, because $\mu_A(I/Q) = d$. We may assume that $i=1$. Let us write $$a_1= \sum_{j=2}^d a_j\xi_j + \sum_{j=1}^d y_j\eta_j$$ with $\xi_j$ and $\eta_j \in A$. Then $\eta_j \in \m$ for all $1 \le j \le d$, since $\sum_{j=1}^d y_j\eta_j \in Q$ and $\mu_A(I/Q) = d$.  Let $c = \sum_{i=1}^d c_i$.
Then $$a_1-\sum_{j=2}^d a_j\xi_j =\sum_{j=1}^d y_j\eta_j \in Q \cap \m^{c- d}=\sum_{j=1}^da_j\m^{c-(d+c_j)}.$$
Hence $$a_1-\sum_{j=2}^d a_j\xi_j = \sum_{j=1}^da_j\rho_j$$ for some $\rho_j \in \m^{c-(d+c_j)}$, so that $a_1(1-\rho_1) \in (a_j \mid 2 \le j \le d)$. Therefore $\rho_1$ is a unit of $A$, since $a_1 \notin (a_j \mid 2 \le j \le d)$. Thus $d=2$ and $c_2 =2$, because $\rho_1 \in \m^{(c_2 + c_3 + \cdots + c_d)-d}$ and $c_j \ge 2$ for all $2 \le j \le d$.

Conversely, assume that $d=2$ and $c_2=2$.
We then have $$I= Q + J = (x_1^{c_1-1},x_1^{c_1-2}x_2, x_2^2).$$
Hence $\mu_A(I)<4=2d$ and so $I \not\subseteq \overline{Q}$ by Proposition \ref{4}.
Thus assertion (1) is proven.
Since $Q \subseteq \m^2$, the second assertion readily follows from Corollary \ref{qi} (3).
\qed
\end{pf}

The following result is the heart of this paper.

\begin{thm}\label{notint}
Let $n = \rmv(A/Q) > 1$ and assume that $I$ is not integral over $Q$.
Then $\e_\m^0(A)\le 2$ and $n=2$.
\end{thm}


\begin{pf}
Firstly, suppose that $d = 1$ and let $a= a_1$. Then $I=(a) + (y_1, y_2, \cdots, y_n)$ and $\m = (a) + (x_1, x_2, \cdots , x_n)$. We have $\mu_A(I) \le n$ by Proposition \ref{4}, because $I$ is not integral over $Q$, while $\mu _{A}(I/Q)=n$ by Lemma \ref{d=0} (1). Hence $I=(y_1, y_2, \cdots, y_n)$ and $a \in \m{\cdot}(y_1, y_2, \cdots, y_n) \subseteq \m^2$. Therefore $\m = (x_1, x_2, \cdots ,x_n)$.
We put $$J:=Q:\m =Q+\m I= Q+(z).$$ Then $\m J = \m Q$ (cf. \cite[Proof of Theorem 2.2]{CP1}; recall that $A$ is not a discrete valuation ring, because $n > 1$). Hence $\mu_A(J) = 2$, because $\ell_A(J/\m J) = \ell_A(J/Q) + \ell_A(Q/\m Q) = 2$. We have $J = \m I = (x_iy_j \mid 1 \le i,j \le n)$, because $Q \subsetneq \m I \subseteq J$.

We divide the proof into two cases.

{\bf Case 1.} ($x_iy_j\notin\m Q$ for some $1 \le i, j \le n$ such that $i \ne j$.)\\
Without loss of generality we may assume that $i=1$ and $j=2$.
Then, because $x_1y_2 \in Q$ but $x_1y_2 \notin \m Q$, we have $Q=(x_1y_2)$. Hence $J=(x_1y_1) + Q  = (x_1y_1, x_1y_2) = x_1{\cdot}(y_1, y_2) \subseteq (x_1)$ because $z \equiv x_1y_1~~ \mod~~ Q$, whence $x_1$ is a non-zerodivisor in $A$. We have $x_1 y_\ell \in \m I = J = x_1{\cdot}(y_1, y_2)$, so that $y_\ell \in (y_1, y_2)$ for all $1 \le \ell \le n$. Thus $I=(y_1, y_2)$. Hence $n = 2$. Because $\m I =  x_1 I$ and $\mu_A(I) = 2$, we have $\m^2 = x_1 \m$, just thanks to the determinantal trick (cf. \cite[Proposition 5.1]{DGH}). Hence $\rme_\m^0(A) = 2$, because $A$ is a Gorenstein local ring of maximal embedding dimension.

{\bf Case 2.} ($x_iy_j \in \m Q$ for all $1 \le i,j \le n$ such that $i \ne j$.)\\
In this case, we have $J = (x_iy_i \mid 1\le i\le n)$, because $J = \m I=(x_iy_j \mid 1\le i,j \le n)$ and $\m J = \m Q$.
Since $\mu_A(J) = 2$, without loss of generality, we may assume that $J=(x_1y_1,x_2y_2)$. Because $x_1y_1-x_2y_2 \notin \m J = \m Q$ and $x_1y_1 \equiv x_2y_2 \equiv z~~ \mod~~ Q$, we have $x_1y_1 = x_2y_2 + a \varepsilon$ with a unit $\varepsilon$ in $A$, while $x_1 y_2 = a\alpha$ and $x_2y_1 = a\beta$ with $\alpha, \beta \in \m$. Hence  $$(x_1+x_2)(y_1-y_2)=a(\varepsilon - \alpha +\beta)$$
with $\varepsilon -\alpha +\beta$ a unit of $A$.
We put  
$$
X_i=
\begin{cases}
x_1+x_2 & (i=1)\\
x_i & (i\ne 1)
\end{cases}
\qquad\text{and}\qquad Y_i=
\begin{cases}
y_1-y_2 & (i=2)\\
y_i & (i\ne 2).
\end{cases}
$$
Then $\m =(X_1, X_2, \cdots,X_n)$, $I=(Y_1, Y_2, \cdots , Y_n)$, and $X_1Y_2 \notin \m Q$ clearly.
Thus thanks to Case 1, we have $n= \rme_\m^0(A)=2$.

Now assume that $d \ge 2$.
Then, by Proposition \ref{4}, we have $\mu_A(I)<n+d$. Since $\mu_A(I/Q) = n$, we may assume that $I=(a_2, a_3, \cdots ,a_d) + (y_1, y_2, \cdots , y_n)$. 
Let $L=(a_2, a_3, \cdots, a_d)$, $\overline A=A/L$, $\overline\m =\m /L$, $\overline Q=Q/L$, and $\overline I=I/L$.
Then $\overline I=\overline Q:\overline\m ^2$ and $\overline A$ is a Gorenstein local ring of dimension $1$ with $\rmv(\overline A/\overline Q)=\rmv(A/Q)=n > 1$.
We have $\mu _{\overline A}(\overline I)\le n$, whence by Proposition \ref{4}, $\overline I$ is not integral over $\overline Q$. Therefore, thanks to the result of the case where $d=1$, we have $n = \rme_{\overline{\m}}^0(\overline{A}) = 2$. We see $\rme_\m^0(A) \le 2$ because $\rme_{\overline{\m}}^0(\overline{A}) \ge \rme_\m^0(A)$, which completes the proof of Theorem \ref{notint}.
\qed
\end{pf}

The following assertion readily follows from Theorem \ref{notint}.

\begin{cor}\label{either}
Suppose that $\e_\m^0(A)\ge 3$. Then $I$ is integral over $Q$, if $n = \rmv(A/Q) > 1$.
\end{cor}

\section{Proof of Theorem 1.1}

Throughout this section let $(A, \m)$ be a Gorenstein local ring with $d=\dim A>0$ and $Q = (a_1, a_2, \cdots ,a_d)$ a parameter ideal in $A$. We put $I=Q:\m^2$. 

The purpose of this section is to prove Theorem 1.1. Let us begin with the following.

\begin{thm}\label{ess}
Suppose that $n = \rmv(A/Q) > 1$ and $I$ is integral over $Q$.
Then 
\begin{enumerate}
\item[$(1)$] $I^3=QI^2$.
\item[$(2)$] $\Gr(I)$ and $\FC(I)$ are Cohen-Macaulay rings.
\end{enumerate}
Hence $\R (I)$ is also a Cohen-Macaulay ring, if $d \ge 3$.
\end{thm}

\begin{pf}
The last assertion directly follows from assertions (1) and (2), because the $a$-invariant $\rma(\Gr(I))$ of $\Gr(I)$ is at most $2-d$ (cf. \cite[THEOREM (1.1), REMARK (3.10)]{GS}).

We may assume that $I^2 \not\subseteq Q$, thanks to Corollary \ref{qi} (1). Choose the element $z\in \m I$ so that $z \in I^2$. Hence $Q: \m = Q+I^2 = Q+(z)$ and so $I^2 = QI + (z)$, because $Q \cap I^2 = QI$ by Corollary \ref{qi} (1). 
Thus $I^3=QI^2+zI$
and we get the required equality $I^3 = QI^2$ modulo the following claim, because $$(Q^2 + zQ) \cap I^3 = (Q^2 \cap I^3) + zQ = Q^2I + zQ \subseteq QI^2$$
 by Corollary \ref{qi} (1). 

\setcounter{claim}{0}
\begin{claim}
$zI \subseteq Q^2 + zQ$.
\end{claim}

\begin{cpf1}
Since $ I = Q + (y_1, y_2, \cdots, y_n)$, it suffices to show that $zy_\ell \in Q^2 + zQ$ for all integers $1 \leq \ell \leq n$. Let $1 \le i \le n$ be an integer such that $i \ne \ell$ and write $z = x_i y_i + q_i$ with $q_i \in \m Q$.
Then $zy_\ell = (x_iy_\ell)y_i + y_\ell q_i \in (\m I)Q$. Since $\m I \subseteq Q:\m = Q+(z)$, we have $zy_\ell \in [Q+(z)]{\cdot} Q = Q^2 + zQ$. Thus $zI \subseteq Q^2 + zQ$.
\qed
\end{cpf1}

As $I^3=QI^2$ and $Q \cap I^2 = QI$ by Corollary \ref{qi} (1), we have $Q \cap I^{i+1}=QI^i$ for every $i \in \Bbb Z$, whence $\Gr(I)$ is a Cohen-Macaulay ring. To show that $\FC(I)$ is a Cohen-Macaulay ring, we need the following.
The equality $\m I^2=\m QI$ in Claim 2 yields, since $I^3=QI^2$, that  the elements  $a_1T ,a_2T, \cdots , a_dT \in \calR(I)$ constitute a regular sequence in $\FC(I)$.

\begin{claim}
$\m I^2=\m QI$.
\end{claim}

\begin{cpf2}
Let $J=(y_1, y_2, \cdots, y_n)$.
Hence $I^2 = QI + J^2$ because $I=Q+J$. It suffices to show that $\m J^2 \subseteq \m QI$. Since $\m = Q+(x_1, x_2, \cdots, x_n)$ and $QJ^2 \subseteq \m QI$, we have only to show $x_\ell y_i y_j \in \m QI$ for all integers $1 \le \ell, i,j \le n$. 
Let us write $x_\ell y_i=\delta _{\ell i}z+q_{\ell i}$ with $q_{\ell i}\in \m Q$. Then $$x_\ell y_iy_j=(\delta_{\ell i}z+q_{\ell i})y_j=\delta_{\ell i}y_j z+q_{\ell i}y_j\in I^3+\m QI =\m QI,$$ because $I^3=QI^2$. Hence $\m I^2=\m QI$.
\qed
\end{cpf2}
\end{pf}

We are now in a position to prove Theorem 1.1.

\begin{tpf}
By Proposition \ref{4}, Corollary 2.8, and Theorem 3.1 we may assume that $n = \rmv(A/Q) = 1$. Hence $\rmv(A)=d+1$.
Let $\m =Q+(x)$ with $x \in \m$; hence $a_1, a_2, \cdots ,a_d, x$ is a minimal basis of $\m$.
We put $$\overline A=A/Q,~\overline\m =\m /Q=(\overline x),~\overline I=I/Q,~\textup{and}~~\ell=\ell_A(\overline A),$$ where $\overline{x} = x~\mod~Q$ be the image of $x$ in $\overline{A}$. Then, since $\overline{\m} = (\overline{x})$, we have $$\ell -1 = \max\{t \in \Bbb{Z} \mid \overline{\m}^t \ne (0)\} \ \ \textup{and}\ \  x^{\ell} \in Q.$$ Hence $\overline{I} = (0) : \overline{\m}^2 =  \overline{\m}^{\ell -2}$ so that $I = Q+\m^{\ell -2} = Q+ (x^{\ell -2})$. Notice that $\ell=\e_Q^0(A) \ge \e_\m^0(A) \ge 3,$ where $\e_Q^0(A)$ denotes the multiplicity of $A$ with respect to $Q$. We then have  $$\m^2 I = [Q\m + (x^2)]{\cdot} [Q+ (x^{\ell -2})] \subseteq \m^2 Q + (x^\ell),$$ because $\m^2 = Q\m + (x^2)$ and $\ell \ge 3$.  Consequently, in order to see that $\m^2 I = \m^2 Q$, it suffices to show the following.

\setcounter{claim}{0}
\begin{claim*}
$x^\ell \in \m^2 Q$.
\end{claim*}

\begin{cpf}
Let us write $x^\ell = \sum_{i=1}^{d} a_i w_i$ with $w_i \in A$.
Let $\widehat A$ be the $\m$-adic completion of $A$ and take an epimorphism $\varphi : B \to \widehat{A}$, where $(B,\n)$ is a regular local ring of dimension $d + 1$. Then $\operatorname{Ker} \varphi$ is a principal ideal in $B$  generated by a single element $\xi \in \n^e$ such that $\xi \notin \n^{e+1}$ where $e = \e_\m^0(A)$; hence $\operatorname{Ker} \varphi \subseteq \n^e$.
Choose elements $\{A_i\}_{1 \le i \le d}, X,$ and $\{W_i\}_{1 \le i \le d}$ of $B$ such that they are the preimages of $\{a_i\}_{1 \le i \le d}, x,$ and $\{w_i\}_{1 \le i \le d}$, respectively.
Then we have $\n =(A_1, A_2, \cdots,A_d, X)$ and $X^\ell - \sum_{i=1}^{d} A_i W_i \in \operatorname{Ker} \varphi \subseteq \n^e$. Hence $\sum_{i=1}^{d} A_i W_i \in \n^e$, because $\ell \geq e$. Consequently, since $(A_1, A_2, \cdots,A_d) \cap \n^e = (A_1, A_2, \cdots,A_d){\cdot}\n^{e-1}$, we see that  $\sum_{i=1}^{d} A_i W_i = \sum_{i=1}^{d} A_i V_i$ for some elements $V_i \in \n^{e-1}$, whence $x^\ell = \sum_{i=1}^{d} a_i v_i$ where $v_i = \varphi (V_i)$. Thus $x^\ell \in Q\m^{e-1} \subseteq Q\m^2$ as is wanted, because $e \ge 3$.
\qed
\end{cpf}

Since $\m^2 I = \m^2 Q$, we have $Q \cap I^2=QI$ similarly as in the proof of Corollary \ref{qi} (1). Therefore, to finish the proof of Theorem 1.1,  we may assume $I^2 \not\subseteq Q$. Since $x^{\ell} \in Q$ and $I^2 = QI + (x^{2\ell - 4})$, we have $2\ell - 4 < \ell $ whence $\ell = e= 3$, so that $I = Q + (x) = \m$. Thus $\m^3 = \m^2 I = Q \m^2$ and so $\Gr(\m) = \FC(\m)$ is a Cohen-Macaulay ring. As $\rma(\Gr(\m)) \le 2-d$, $\calR(\m)$ is a Cohen-Macaulay ring if $d \ge 3$. This completes the proof of Theorem 1.1.
\qed
\end{tpf}

\section{Examples}

In this section we explore three examples to show that the non-Gorenstein case is rather wild.

\begin{ex}
Let $n \ge 2$ be an integer and let 
$$A = k[[X_1,X_2,\cdots,X_n]] / (X_iX_j \mid 1\le i < j \le n)$$
 where $ k[[X_1,X_2,\cdots,X_n]]$ denotes the formal power series ring over a field $k$.
Then $A$ is a one-dimensional reduced  local ring with $\e_\m^0(A) = n$. For every parameter ideal $Q$ in $A$, we have 
$$Q:\m^2 \not\subseteq \overline{Q},$$
where $\overline{Q}$ denotes the integral closure of $Q$.
\end{ex}
\begin{pf}
Let $I=Q:\m^2$ and assume that $I \subseteq \overline{Q}$. We write $Q=(a)$.
Then $a= \sum_{i=1}^n x_i^{c_i}\varepsilon_i$ for some units $\varepsilon_i$ in $A$ and some integers $c_i \ge 1$.
Let $1\le i\le n$ be an integer. If $c_i\ge 2$, we then have $x_i^{c_i-1} \in I$ but $x_i^{c_i-1}$ is not integral over $Q$.
Hence $c_i = 1$ for all $1 \le i\le n$ and so $a =\sum_{i=1}^n x_i\varepsilon_i$. Therefore $\m^2 = Q\m$ so that we have $I = A$, which is absurd.
\qed
\end{pf}

Letting $n = 2$, this Example 4.1 shows the assumption that $\e_\m^0(A) \ge 3$ in Theorem 1.1 is not superfluous.

It seems natural and quite interesting to ask what happens in the case where $A$ is a numerical semi-group ring. Let us explore one example.

\begin{ex}
Let $H=\left<4,7,9\right>$ be the numerical semi-group generated by $4, 7$, and $9$ and let $A = k[[t^4,t^7,t^9]]~\subseteq~k[[t]]$, where $V= k[[t]]$ denotes the formal power series ring with one indeterminate $t$ over a field $k$. Then $A$ is a one-dimensional non-Gorenstein Cohen-Macaulay local ring.
Let $0 < s \in H$. We put $Q=(t^s)$ and  $I=Q:\fkm^2$. 
Then $I \subseteq \overline{Q}$. We have $I \subseteq \m^2$ if $s \ge 11$, whence $I^2 \subseteq Q$. However
\begin{enumerate}
\item[$(1)$] $\rmr_Q(I) = 
\left\{
\begin
{array}{lll}
1 & \text{ if } & s=9,\\
2 & \text{ if } & s=4,~8,~\textup{or}~s \ge 11,\\
3 & \text{ if } & s=7.
\end{array}
\right.$
\item[$(2)$] $\rmG(I)$ is a Cohen-Macaulay ring if and only if $s=4,8,9$.
\item[$(3)$] $\rmF(I)$ is a Cohen-Macaulay ring if and only if $s=4, 9$.
\item[$(4)$] $\rmF(I)$ is always a Buchsbaum ring.
\item[$(5)$] $\rmG(I)$ is a Buchsbaum ring if and only if $s\ne 7$.
\item[$(6)$] $\m^2 I \ne \m^2 Q$ if $s =8,  11$.
\end{enumerate}
\end{ex}

\begin{pf}
We have $n \in H$ for all integers $n \ge 11$ but $10 \not\in H$. Hence the conductor of $H$ is 11. Notice that $t^n \in \m^2$ for all $n \in \Bbb Z$ such that $n \ge 11$, where $\m = (t^4, t^7, t^9)$ denotes the maximal ideal in $A$. Hence $I \subseteq \overline Q$. In fact, let $n \in H$ and assume that $t^n \in I$ but $n < s$. Then $t^{s-n+10} \in \m^2$ because $s-n+10 \ge 11$, so that $t^{s+10} = t^n t^{s-n+10} \in Q = (t^s)$ whence $t^{10} \in A$, which is impossible. Thus, for every $n \in H$ with $t^n \in I$, we have $t^n \in t^sV \cap A = \overline Q$, whence $I \subseteq \overline Q$ (recall that  $I$ is a monomial ideal generated by the elements $\{t^n \mid n \in H~\textup{such that}~t^n \in I\}$). In particular we have  $I \subseteq \m^2$ if $s \ge 11$, whence $I^2 \subseteq Q$.

We note the following.
\setcounter{claim}{0}
\begin{claim}
Let $s_2 \ge s_1 \ge 11$ be integers and let $q = s_2 -s_1$. We put $Q_i = (t^{s_i})$ and $I_i = Q_i : \m^2$ for $i = 1,2$. Then we have the following. 
\begin{enumerate}
\item[$(1)$] $I_2 = t^{q}I_1$. 
\item[$(2)$] $\calR(I_1) \cong \calR(I_2)$ as graded $A$-algebras. 
\item[$(3)$] $\rmF(I_1) \cong \rmF(I_2)$ as graded $A/\m$-algebras.
\item[$(4)$] $\rmr_{Q_1}(I_1) = \rmr_{Q_2}(I_2)$.
\end{enumerate}
\end{claim}

\begin{cpf1}
Let $\varphi = \widehat{t^q} : V \to V$ be the $V$-linear map defined by $\varphi (x ) = t^qx$ for all $x \in V$. Then, since $\varphi (Q_1) = Q_2$ and $\varphi (I_1) \subseteq I_2$, the map $\varphi$ induces a monomorphism $$\xi : I_1/Q_1 \to I_2/Q_2, ~~\ \ x~\mod~Q_1~\mapsto~t^qx~\mod~Q_2$$ of $A$-modules. As $I_i/Q_i \cong \operatorname{Ext}_A^1(A/\m^2, A)$ (recall that $t^{s_i} \in \m^2$), we see $\ell_A(I_1/Q_1) = \ell_A(I_2/Q_2)$, whence $\xi : I_1/Q_1 \to I_2/Q_2$ is an isomorphism, so that $\varphi(I_1) = I_2$. Thus assertion (1) follows. Notice that $$\calR(I_2) = A[(t^qI_1){\cdot}T] =A[I_1{\cdot}t^qT]~~\textup{and}~~\calR(I_1) = A[I_1T]$$ with $T$ an indeterminate over $A$. Then, since $t^qT$ is also transcendental over the ring $A$, we get an isomorphism $\xi : \calR(I_1) \to  \calR(I_2)$ of graded $A$-algebras such that  $\xi(t^{s_1}T)=t^{s_2}T$. Hence we have assertion (2). Because $\rmF(I_i) = \calR(I_i)/\m \calR(I_i)$, we  readily have an isomorphism $\eta : \rmF(I_1) \to \rmF(I_2)$ of graded $A/\m$-algebras such that $\eta (\overline{t^{s_1}T}) = \overline{t^{s_2}T}$, where $\overline{t^{s_i}T}$ denotes the image of $t^{s_i}T$ in $\rmF(I_i)$. Hence assertion (4) also follows, because $$\rmr_{Q_i}(I_i) = \max \{ n \in \Bbb Z \mid [\rmF(I_i)/(\overline{t^{s_i}T})]_n \ne (0) \}, $$
where $[\rmF(I_i)/(\overline{t^{s_i}T})]_n$ denotes the homogeneous component of the graded ring $\rmF(I_i)/(\overline{t^{s_i}T})$ of degree $n$.
\qed
\end{cpf1}

We put $\calR = \calR(I)$, $G = \rmG(I)$, and $F = \rmF(I)$. Let $M = \m \calR + \calR_+$ be  the graded maximal ideal in $\calR$ and we denote by $\rmH_M^0(*)$ the $0^{\underline{th}}$ local cohomology functor with respect to $M$. Let $a =t^s$ and $f = aT \in \calR=A[IT]$. For each graded $\calR$-module $L$, let $[\rmH_M^0(L)]_n$ $(n \in \Bbb Z)$ denote the homogeneous component of $\rmH_M^0(L)$ of degree $n$. 

Let $\widetilde{I} =\bigcup_{n \ge 0}[I^{n+1} : I^n]$ denote the Ratliff-Rush closure of $I$. The following assertions readily follow from the equalities  that $$\widetilde{I} = \bigcup_{n \ge 0}[I^{n+1} : a^n]~~\ \textup{and}~~\ I^{n + \ell} =a^\ell I^n$$ for all integers $n \ge r=\rmr_Q(I)$ and $\ell \ge 1$, whose details are  left to the reader.
\begin{claim}
Let  $r = \rmr_Q(I)$. Then 
\begin{enumerate} 
\item[$(1)$] $[\rmH_M^0(G)]_n =(0)$ for all $n \ge r-1$.
\item[$(2)$] $[\rmH_M^0(F)]_n =(0)$ for all $n \ge r$.
\item[$(3)$] Suppose that $r = 2$.  Then $\widetilde{I} = I^2 : a$ and $[\rmH_M^0(G)]_0 \cong \widetilde{I}/I$ as $A$-modules.
\end{enumerate}
\end{claim}

We now consider the case $s = 11$. We then have $I = (t^{11}, t^{12}, t^{14}, t^{17})$, $I^3 = QI^2$, and $$I^2  = QI + (t^{24}) \ne QI$$ since $t^{24} \not\in QI = (t^{22}, t^{23}, t^{25}, t^{28})$. Hence $\rmr_Q(I) = 2$. Because $\widetilde{I} = I : t^{11} = I + (t^{13}) \ne I$ and $$\rmH_M^0(G) = [\rmH_M^0(G)]_0 \cong \widetilde{I}/I$$
by Claim 2 (3),  we see that $G$ is not a Cohen-Macaulay ring but a Buchsbaum ring with $\ell_A(\rmH_M^0(G)) = \ell_A(\widetilde{I}/I) = 1$.
Notice that $\m^2 I = (t^{19}, t^{20}, t^{22}, t^{25}) \ne \m^2 Q = (t^{19}, t^{22}, t^{24}, t^{25})$. Since $t^{11}t^{17}=t^{28} =t^{4}t^{24}\in \m I^2$ but $t^{17} \not\in \m I = (t^{15}, t^{16},t^{18},t^{21})$, the element $f = t^{11}T \in \calR$ is a zerodivisor in $F$, whence $F$ is a Buchsbaum ring by Claim 2 (2) but not a Cohen-Macaulay ring.

If $s > 11$, then thanks to Claim 1 and the assertions in the case where $ s = 11$, we have $\rmr_Q(I) = 2$ and $F$ is a Buchsbaum ring but not Cohen-Macaulay. To see that $G$ is a Buchsbaum ring, recall that $$\rmH_M^0(G) = [\rmH_M^0(G)]_0 \cong \widetilde{I}/I$$ since $\rmr_Q(I) = 2$. Let $Q' = (t^{11})$ and $I' = Q' : \m^2$. Then because $\widetilde{I'} = {I'}^2 : t^{11}$ and $\widetilde{I} = I^2 : t^{s}$ (see Claim 1 (4)), it is standard to check that $t^{s - 11}{\cdot}\widetilde{I'} = \widetilde{I}$, so that we have $\ell_A(\widetilde{I}/I) = \ell_A(\widetilde{I'}/I')= 1$ (recall that $t^{s - 11}{\cdot}I' = I$; cf. Claim 1 (1)), whence $G$ is a Buchsbaum ring with $\ell_A(\rmH_M^0(G))=1$.

Let $s = 4$. Then $I =\m$. The ring $G~(=F)$ is a Cohen-Macaulay ring, since $\m^3 = Q\m^2$ and $Q \cap \m^2 = Q\m$.

Let $s = 7$, then $I = (t^{7}, t^{8}, t^{11}, t^{13})$, $I^2 =(t^{14}, t^{15}, t^{16}) \subseteq Q$, and $t^{16} \notin QI = (t^{14},t^{15},t^{20})$. Hence $G$ is not a Cohen-Macaulay ring. We have $I^4 = QI^3$ but $I^3 = QI^2 + (t^{24}) \ne QI^2$. Hence $\rmr_Q(I) = 3$. Because $t^{7}t^{13}=t^{20} =t^4t^{16}\in \m I^2$ but $t^{13} \notin \m I =(t^{11}, t^{12}, t^{14}, t^{17})$, $F$ is not a Cohen-Macaulay ring. To see that $G$ is not a Buchsbaum ring, let $W = \rmH_M^0(G)$. Then $W = W_0 + W_1$ by Claim 2 (1). It is now direct to check that $W_0 = \{\overline{c} \mid c \in (t^9)\}$ and $W_1 = \{\overline{cT} \mid c \in (t^{17}) \}$ where $\overline{*}$ denotes the image of the corresponding element of $\calR$
in $G$. Because $\overline{t^9} \ne 0$ in $G$ and $t^9{\cdot}I=(t^{16}, t^{17}, t^{20}, t^{22}) \not\subseteq I^2 = (t^{14}, t^{15}, t^{16})$, we see $MW_0 \ne (0)$, whence $G$ is not a Buchsbaum ring. Similarly, one can directly check that $$\rmH_M^0(F) = [\rmH_M^0(F) ]_1 =\{\overline{cT} \mid c \in (t^{13})\} \cong A/\m ,$$ so that $F$ is a Buchsbaum ring but not Cohen-Macaulay. 

Let $s=8$. Then $I=(t^{8}, t^{9}, t^{11}, t^{14})$ and  $I^3 = QI^2$. We have $\m^2 I = (t^{16},t^{17},t^{19},t^{22}) \ne \m^2 Q = (t^{16}, t^{19}, t^{21}, t^{22})$ and $\ell_A(\m^2 I /\m^2 Q) = 1$. To see that $G$ is a Cohen-Macaulay ring, we have only to show that $Q \cap I^2 = QI$. Since $I^2 = QI + (t^{18})$, we have $Q\cap I^2 = QI + [Q \cap (t^{18})]$.@Let $\varphi \in Q \cap (t^{18})$ and write $\varphi = t^8\xi = t^{18}\eta$ with $\xi, \eta \in A$. Then $\xi = t^{10}\eta$. Because $10 \notin H=\left<4, 7, 9\right>$, we have $\eta \in \m$ so that $\varphi = t^{18}\eta \in t^{18}\m = (t^{22}, t^{25}, t^{27}) \subseteq QI = (t^{16}, t^{17}, t^{19}, t^{22})$. Hence $Q \cap I^2 = QI$ and $G$ is a Cohen-Macaulay ring.
The ring $F$ is Buchsbaum by Claim 2 (2)  but not a Cohen-Macaulay ring, because $t^8t^{14} = t^{22}= t^4(t^{9})^2 \in \m I^2$ but $t^{14} \notin \m I =(t^{12}, t^{13}, t^{15}, t^{18})$.

Let $s = 9$. Then $I = (t^{9}, t^{12}, t^{14}, t^{15})$ and $I^2 = QI$, whence $G$ and $F$ are both Cohen-Macaulay rings. This completes the proofs of all the assertions. \qed
\end{pf}

Our last example shows that unless $A$ is Gorenstein, the reduction number $\rmr_Q(I)$ can be arbitrarily large even if $I \subseteq \overline{Q}$, where $I=Q:\m^2$ and $\overline{Q}$ denotes the integral closure of $Q$.

\begin{ex}
Let $n \ge 3$ be an integer and let 
$$
a_i=
\begin{cases}
2n-1 & \quad (i=1), \\
(2n+1)i-2n-2 & \quad (2\le i\le n).
\end{cases}
$$
Let $H = \left< a_1, a_2, \cdots , a_n \right>$ be the numerical semi-group generated by $a_i$'s. Let $A= k[[t^{a_1},t^{a_2}, \cdots,t^{a_n}]]  \subseteq k[[t]]$ be the semi-group ring of $H$, where $k[[t]]$ denotes the formal power series ring with one indeterminate $t$ over a field $k$. Then $A$ is a one-dimensional Cohen-Macaulay local ring with the maximal ideal $\m =(t^{a_1},t^{a_2}, \cdots,t^{a_n})$. Let $Q=(t^{2a_1})$ and $I=Q:\m^2$. Then $I \subseteq \overline{Q}$ and $\rmr_Q(I) = 2n-2$.
\end{ex}

\begin{pf}
Let $B=k[[X_1,X_2, \cdots,X_n]]$ $(n \ge 2)$ be the formal power series ring over the field $k$ and let 
$$\varphi : B \to  k[[t^{a_1},t^{a_2}, \cdots,t^{a_n}]]$$ be the homomorphism of $k$-algebras defined by $\varphi(X_i) = t^{a_i}$ for all $1 \le i \le n$. 
Let $\rmI_2(M)$ be the ideal in $B$ generated by all the $2 \times 2$ minors of the following matrix 
$$
M=
\begin{pmatrix}
X_1 & X_2 & X_3 & \cdots & X_{n-1} & X_n \\
X_2^2 & X_3 & X_4 & \cdots & X_n & X_1^{n+1}
\end{pmatrix}.
$$
We then have $\operatorname{Ker} \varphi = \rmI_2(M)$, because $\ell_B(B/[\rmI_2(M) + (X_1)]) = 2n-1 = a_1$. Let us identify  $A=B/\rmI_2(M) = k[[t^{a_1},t^{a_2}, \cdots,t^{a_n}]]$. Let  $x_i = X_i~\mod~ \rmI_2(M)$ be the image of $X_i$ in $B/\rmI_2(M)$ for each $1\le i \le n$; hence $\m=(x_1,x_2, \cdots,x_n)$. With this notation it is standard and easy to check that $\m^2=(x_1,x_2)\m$, $I=x_1\m +(x_2x_n),$  and $I^i=x_1^i\m^i$ for all $i\ge 2$.

Recall now that $\m^{2n-1} = x_1\m^{2n-2}$, because  $(x_1)$ is a minimal reduction of $\m$ and $\e_\m^0(A) = 2n-1$.
Hence $$I^{2n-1} = x_1^{2n-1}\m^{2n-1} = x_1^{2n-1}{\cdot}x_1\m^{2n-2} = x_1^2{\cdot}x_1^{2n-2}\m^{2n-2} = x_1^2 I^{2n-2} = QI^{2n-2}.$$ We must show that $I^{2n-2}\ne QI^{2n-3}$. To see this, we explore the following system of generators of $\m^{2n-3}$;
\begin{align*}
\m^{2n-3} & =(t^{a_1},t^{a_2},\cdots,t^{a_n})^{2n-3} \\
& = (t^{\sum_{i=1}^n c_i a_i} \mid c_i \ge 0 ,~\sum_{i=1}^n c_i = 2n-3)\\
& = (t^{(2n-i-3)a_1+ia_2} \mid 0\le i \le 2n-3)\\
&\ \ \ \  + (t^{\sum_{i=1}^n c_i a_i} \mid c_j > 0~\textup{for some}~j \ge 3,~\sum_{i=1}^n c_i = 2n-3).
\end{align*}
Notice that $\{(2n-i-3)a_1+ia_2 = (4n^2 -8n + 3) + i\}_{0\le i \le 2n-3}$ are continuous integers and that $$\sum_{i=1}^n c_i a_i \ge (2n-4)a_1 + a_3  = 4n^2 -6n +5,$$ if $c_j > 0$ for some $j \ge 3$ and $\sum_{i=1}^n c_i = 2n-3$. Hence $$\m^{2n-3} \subseteq (t^i \mid 4n^2 -8n + 3 \le i \le 4n^2-6n) + (t^i \mid i\in H,~i \ge 4n^2 - 6n + 5).$$ Therefore
\begin{align*}
(\sharp) \ \ \ t^{2n-1}\m^{2n-3} & \subseteq (t^i \mid 4n^2-6n+2\le i\le 4n^2-4n-1) +(t^i \mid i\in H,~i \ge 4n^2 - 4n + 4).
\end{align*}
Suppose now that $I^{2n-2} = QI^{2n-3}$. Then $\m^{2n-2} = x_1\m^{2n-3}$, since $n \geq 3$ (recall that $I^i=x_1^i\m^i$ for all $i\ge 2$); hence $x_2^{2n-2} \in x_1\m^{2n-3}$.
Recall that $x_1 = t^{2n-1}$ and $x_2 = t^{2n}$. Then, because  $4n^2-4n-1<4n^2-4n<4n^2-4n+4$, we get by $(\sharp)$ that $$t^{4n^2-4n} \in \m{\cdot}(t^i \mid 4n^2-6n+2\le i\le 4n^2-4n-1),$$
which  is however impossible, since $$\m{\cdot}(t^i \mid 4n^2-6n+2\le i\le 4n^2-4n-1) \subseteq t^{4n^2-4n+1}k[[t]]$$
(recall that $a_i+(4n^2-6n+2) \ge a_1+(4n^2-6n+2) = 4n^2-4n+1$ for all $1 \le i \le n$).
This is the required contradiction and we conclude that $I^{2n-2}\ne QI^{2n-3}$. Thus $\rmr_Q(I) = 2n-2$.
\qed
\end{pf}



\begin{thebibliography}{GSa3}


\bibitem[CHV]{CHV} A. Corso, C. Huneke, and W. V. Vasconcelos,
{\it On the integral closure of ideals},
Manuscripta Math., {\bf 95} (1998), 331--347.


\bibitem[CP1]{CP1} A. Corso and C. Polini,
{\it Links of prime ideals and their Rees algebras,}
J. Algebra, {\bf 178} (1995), 224--238.

\bibitem[CP2]{CP2} A. Corso and C. Polini, 
{\it Reduction number of links of irreducible varieties,}
J. Pure Appl. Algebra, {\bf 121} (1997), 29--43.

\bibitem[CPV]{CPV} A. Corso,  C. Polini, and W. V. Vasconcelos,
{\it Links of prime ideals,}
Math. Proc. Camb. Phil. Soc., {\bf 115} (1994), 431--436.


\bibitem[DGH]{DGH} M. D'anna, A. Guerrieri, and W. Heinzer, {\it Ideals having a one-dimensional fiber cone}, Ideal Theoretic Methods in Commutative Algebra, 155--170,  Lecture Notes in Pure and Appl. Math., 220, Dekker, New York, 2001.

\bibitem[G]{G} S. Goto,
{\it Integral closedness of complete intersection ideals,}
J. Algebra, {\bf 108}(1987), 151--160.

\bibitem[GN]{GN} S. Goto and K. Nishida, {\it Hilbert coefficients and Buchsbaumness of associated graded rings}, J. Pure Appl. Algebra, {\bf 181} (2003), 61--74.

\bibitem[GSa1]{GSa1} S. Goto and H. Sakurai, {\it The equality $I\sp 2=QI$ in Buchsbaum rings}, Rend. Sem. Mat. Univ. Padova, {\bf  110} (2003), 25--56.

\bibitem[GSa2]{GSa2} S. Goto and H. Sakurai, {\it The reduction exponent of socle ideals associated to parameter ideals in a Buchsbaum local ring of multiplicity two},  J. Math. Soc. Japan, {\bf 56} (2004), 1157--1168.

\bibitem[GSa3]{GSa3} S. Goto and H. Sakurai, {\it When does the equality $I\sp 2=QI$ hold true in Buchsbaum rings?},  Commutative Algebra, 115--139, Lect. Notes Pure Appl. Math., {\bf 244}, 2006.
  

\bibitem[GS]{GS} S. Goto and Y. Shimoda, {\it On the Rees algebras of Cohen-Macaulay local rings}, Commutative Algebra (Fairfax, Va., 1979), 201--231, Lecture Notes in Pure and Appl. Math., 68, Dekker, New York, 1982.

\bibitem[NR]{NR} D. G. Northcott and D. Rees,  {\it Reductions of ideals in local rings},  Proc. Camb. Phil. Soc.,  {\bf 50} (1954), 145--158.

\bibitem[PU]{PU} C. Polini and B. Ulrich, 
{\it Linkage and reduction numbers,}
Math. Ann., {\bf 310} (1998), 631-651.

\bibitem[W]{W} H.-J. Wang, 
{\it Links of symbolic powers of prime ideals,}
Math. Z., {\bf 256} (2007), 749-756.
\end{thebibliography}
\end{document}